\documentclass{amsart}
\usepackage{amssymb}

\newtheorem{theorem}{Theorem}
\theoremstyle{plain}

\newtheorem{proposition}{Proposition}

\numberwithin{equation}{section}
\input{tcilatex}

\begin{document}
\title[Noncommutative BV-geometry and Matrix integrals]{Noncommmutative
Batalin-Vilkovisky geometry and Matrix integrals.\vspace{-0.35cm}}
\author{Serguei Barannikov\vspace{-0.3cm}}
\address{Ecole Normale Superieure, 45, rue d'Ulm 75230, Paris, France }
\email{sergueibar@gmail.com}

\begin{abstract}
We associate the new type of \emph{supersymmetric} matrix models with any
solution to the quantum master equation of the noncommutative
Batalin-Vilkovisky geometry. The asymptotic expansion of the matrix
integrals gives homology classes in the Kontsevich compactification of the
moduli spaces, which we associated with the solutions to the quantum master
equation in our previous paper. We associate with the \emph{queer} matrix
superalgebra equipped with an odd \emph{differentiation}, whose square is
nonzero, the family of cohomology classes of the compactification. This
family is the generating function for the products of the tautological
classes. The simplest example of the matrix integrals in the case of
dimension zero is a supersymmetric extenstion of the Kontsevich model of
2-dimensional gravity. 
\end{abstract}

\maketitle

\vspace{-0.05in}\emph{Notations}.I work in the tensor category of super
vector spaces, over an algebraically closed field $k$, $char(k)=0$. Let $%
V=V^{even}\oplus V^{odd}$ be a super vector space. We denote by $\overline{%
\alpha }$ the parity of an element $\alpha $ and by $\Pi V$ the super vector
space with inversed parity. For a finite group $G$ acting on a vector space $%
U$, we denote via $U^{G}$ the space of invariants with respect to the action
of $G$. Element $(a_{1}\otimes a_{2}\otimes \ldots \otimes a_{n})$ of $%
A^{\otimes n}$ is denoted by $(a_{1},a_{2},\ldots ,a_{n})$. Cyclic words,
i.e. elements of the subspace $(V^{\otimes n})^{\mathbb{Z}/n\mathbb{Z}}$ are
denoted via $(a_{1}\ldots a_{n})^{c}$

\section{Noncommutative Batalin-Vilkovisky geometry.\protect\footnotetext[1]{%
Preprint NI06043, 25/09/2006, Isaac Newton Institute for Mathematical
Sciences; preprint HAL-00102085.}\protect\footnotetext[2]{%
The paper submited to the "Comptes rendus" of the French Academy of Science
on May,17,2009; presented for publication by Academy member M.Kontsevich on
May,20,2009. }}

\subsection{\protect\bigskip Even inner products.}

Let $B:V^{\otimes 2}\rightarrow k$ be an even symmetric inner product on $V$%
: 
\begin{equation*}
B(x,y)=(-1)^{\overline{x}\overline{y}}B(y,x)
\end{equation*}%
I introduced in \cite{B1} the space $F=\bigoplus_{n=1}^{n=\infty }F_{n}$ 
\begin{equation}
F_{n}=((\Pi V)^{\otimes n}\otimes k[\mathbb{S}_{n}]^{\prime })^{\mathbb{S}%
_{n}}  \label{Feven}
\end{equation}%
where $k\mathbf{[}\mathbb{S}_{n}]^{\prime }$ denotes the super $k\mathbf{-}$%
vector space with the basis indexed by elements $(\sigma ,\rho _{\sigma })$,
where $\sigma \in \mathbb{S}_{n}$ is a permutation with $i_{\sigma }$ cycles 
$\sigma _{\alpha }$ and $\rho _{\sigma }=\sigma _{1}\wedge \ldots \wedge
\sigma _{i_{\sigma }}$, $\rho _{\sigma }\in Det(Cycle(\sigma ))$, $%
Det(Cycle(\sigma ))=Symm^{i_{\sigma }}(k^{0|i_{\sigma }})$, is one of the
generators of the one-dimensional determinant of the set of cycles of $%
\sigma $, i.e. $\rho _{\sigma }$ is an order on the set of cycles defined up
to even reordering, and $(\sigma ,-\rho _{\sigma })=-(\sigma ,\rho _{\sigma
})$. The group $\mathbb{S}_{n}$ acts on $k[\mathbb{S}_{n}]^{\prime }$ by
conjugation. The space $F$ is naturally isomorphic to: 
\begin{equation*}
F=Symm(\oplus _{j=1}^{\infty }\Pi (\Pi V^{\otimes j})^{\mathbb{Z}/j\mathbb{Z}%
})
\end{equation*}%
The space $F$ carries the naturally defined Batalin-Vilkovisky differential $%
\Delta $ (see loc.cit. and references therein). It is the operator of the
second order with respect to the multiplication and is completely determined
by its action on the second power of $\oplus _{j=1}^{\infty }\Pi (\Pi
V^{\otimes j})^{\mathbb{Z}/j\mathbb{Z}}$. If one chooses a basis $\{a_{i}\}$
in $\Pi V$, in which the antisymmetric even inner product defined by $B$ on $%
\Pi V$ has the form $(-1)^{\overline{a}_{i}}B(\Pi a_{i},\Pi a_{j})=b_{ij}$ ,
then the operator $\Delta $ sends a product of two cyclic words $(a_{\rho
_{1}}\ldots a_{\rho _{r}})^{c}(a_{\tau _{1}}\ldots a_{\tau _{t}})^{c}$, to 
\begin{equation}
\sum_{p,q}(-1)^{\varepsilon _{1}}b_{\rho _{p}\tau _{q}}(a_{\rho _{1}}\ldots
a_{\rho _{p-1}}a_{\tau _{q+1}}\ldots a_{\tau _{q-1}}a_{\rho _{p+1}}\ldots
a_{\rho _{r}})^{c}+  \label{delltaa}
\end{equation}%
\vspace{-0.1in}\vspace{-0.03in}%
\begin{multline*}
+\sum_{p\pm 1\neq q\func{mod}r}(-1)^{\varepsilon _{2}}b_{\rho _{p}\rho
_{q}}(a_{\rho _{1}}\ldots a_{\rho _{p-1}}a_{\rho _{q+1}}\ldots a_{\rho
_{r}})^{c}(a_{\rho _{p+1}}\ldots a_{\rho _{q-1}})^{c}(a_{\tau _{1}}\ldots
a_{\tau _{t}})^{c} \\
+\sum_{p\pm 1\neq q\func{mod}r}(-1)^{\varepsilon _{3}}b_{\tau _{p}\tau
_{q}}(a_{\rho _{1}}\ldots a_{\rho _{r}})^{c}(a_{\tau _{1}}\ldots a_{\tau
_{p-1}}a_{\tau _{q+1}}\ldots a_{\tau _{t}})^{c}(a_{\tau _{p+1}}\ldots
a_{\tau _{q-1}})^{c}
\end{multline*}%
where $\varepsilon _{i}$ are the standard Koszul signs, which take into the
account that the parity of any cycle is opposite to the sum of parities of $%
a_{i}$ : $\overline{(a_{\rho _{1}}\ldots a_{\rho _{r}})^{c}}=1+\sum 
\overline{a_{\rho _{i}}}$. It follows from the loc.cit., prop. 2, that $%
\Delta $ defines the structure of Batalin-Vilkovisky algebra on $F$, in
particular $\Delta ^{2}=0$. The solutions of \emph{the quantum master
equation} in $F$ 
\begin{equation}
\hbar \Delta S+\frac{1}{2}\{S,S\}=0,\,\,S=\sum_{g\geq 0}\hbar
^{2g-1+i}S_{g,i\,,n\,},\,\,S_{g,i,n}\in Symm^{i}\cap F_{n}^{even},
\label{qms}
\end{equation}%
with $S_{0,1,1}=0$, are in one-to one correspondence, by the loc.cit.,
theorem 1, with the structure of $\mathbb{Z}/2\mathbb{Z}$-graded \emph{%
quantum }$A_{\infty }-$\emph{algebra }on $V$,\emph{\ }i.e.the algebra over
the $\mathbb{Z}/2\mathbb{Z}$-graded modular operad $\mathcal{F}_{\mathcal{K}}%
\mathbb{S}$, where $\mathbb{S}$ is the $\mathbb{Z}/2\mathbb{Z}$-graded
version of the twisted modular $\mathcal{K}-$operad $\widetilde{\mathfrak{s}}%
\Sigma \mathbb{S}[t]$, with components $k[\mathbb{S}_{n}^{\prime }][t]$,
described in loc.cit. The $(g=0,i=1)-$part is the cyclic $\mathbb{Z}/2%
\mathbb{Z}$-graded $A_{\infty }-$algebra with the \emph{even} invariant
inner product on $Hom(V,k)\overset{B}{\simeq }V$ . Recall, see loc.cit.,
that for any solution to (\ref{qms}), with $S_{0,1,1}=S_{0,1,2}=0$, the
value of partition function $c_{S}(G)$ on a stable ribbon graph $G$, with no
legs, is defined by contracting the product of tensors $\bigotimes_{v\in
Vert(G)}S_{g(v),i(v),n(v)}$ with $B^{\otimes Edge(G)}$ with appropriate
signs.

\begin{proposition}
(\cite{B1}, s.10,11) The graph complex $\mathcal{F}_{\mathcal{K}}\mathbb{S}%
((0,\gamma ,\nu )))$ (part of $\mathcal{F}_{\mathcal{K}}\mathbb{S}$ with no
legs) is naturally identified with the cochain CW-complex $C^{\ast }(%
\overline{\mathcal{M}}_{\gamma ,\nu }^{\prime }/\mathbb{S}_{\nu })$ of the
Kontsevich compactification of the moduli spaces of Riemann surfaces from (%
\cite{K}). For any solution to the quantum master equation (\ref{qms}) in $F$%
, with $S_{0,1,1}=S_{0,1,2}=0$, the partition function on stable ribbon
graphs $c_{S}(G)$ defines the characterestic homology class in $H_{\ast }(%
\overline{\mathcal{M}}_{\gamma ,\nu }^{\prime }/\mathbb{S}_{\nu })$.
\end{proposition}

\subsection{Odd inner products}

Let $V=V^{0}\oplus V^{1}$ be a super vector space and $B$ be an odd
symmetric inner product, $B:V^{\otimes 2}\rightarrow \Pi k$, $B(x,y)=(-1)^{%
\overline{x}\overline{y}}B(y,x)$. The analog of the space $F$ in this
situation has components 
\begin{equation*}
\widetilde{F}_{n}=(V^{\otimes n}\otimes k[\mathbb{S}_{n}])^{\mathbb{S}_{n}}
\end{equation*}%
where $k\mathbf{[}\mathbb{S}_{n}]$ is the group algebra of $\mathbb{S}_{n}$,
and $\mathbb{S}_{n}$ acts on $k[\mathbb{S}_{n}]$ by conjugation. The space $%
\widetilde{F}$ is naturally isomorphic in this case to: 
\begin{equation}
\widetilde{F}=Symm(\oplus _{j=1}^{\infty }(V^{\otimes j})^{\mathbb{Z}/j%
\mathbb{Z}})  \label{Fodd}
\end{equation}%
The space $\widetilde{F}$ carries the naturally defined second order
differential defined by the formula (\ref{delltaa}) with $a_{i}\in V$, $%
b_{ij}=B(a_{i},a_{j})$ and the Koszul signs $\varepsilon _{i}$, which now
correspond to the standard parity of cycles :$\overline{(a_{\rho _{1}}\ldots
a_{\rho _{r}})^{c}}=\sum \overline{a_{\rho _{i}}}$ . Again, it follows from
the loc.cit., prop. 2, that $\Delta ^{2}=0$, and that $\Delta $ defines the
structure of Batalin Vilkovisky algebra on $\widetilde{F}$.

The solutions of quantum master equation (\ref{qms}) in $\widetilde{F}$,
with $S_{0,1,1}=0$, are in one-to one correspondence, by the loc.cit., with
the structure of algebra over the twisted modular operad $\mathcal{F}%
\widetilde{\mathbb{S}}$ on the vector space $V$ . Here $\widetilde{\mathbb{S}%
}$ is the\emph{\ untwisted} $\mathbb{Z}/2\mathbb{Z}$-graded version of $%
\widetilde{\mathfrak{s}}\Sigma \mathbb{S}[t]$. The components $\widetilde{%
\mathbb{S}}((n))$ are the spaces $k\mathbf{[}\mathbb{S}_{n}][t]$, with the
composition maps defined as in loc.cit.,sect. 9. The Feynman transform $%
\mathcal{F}\widetilde{\mathbb{S}}$ is a $\mathcal{K}-$twisted modular
operad, whose $(g=i=0)$-part corresponds to the cyclic $A_{\infty }-$algebra
with the \emph{odd} invariant inner product on $Hom(\Pi V,k)\overset{B}{%
\simeq }V$.

\begin{proposition}
(\cite{B1},s.10,11) The graph complexes $\mathcal{F}\widetilde{\mathbb{S}}%
((0,\gamma ,\nu )))$ (part of $\mathcal{F}\widetilde{\mathbb{S}}$ with no
legs) are naturally identified with the cochain CW-complexes $C^{\ast }(%
\overline{\mathcal{M}}_{\gamma ,\nu }^{\prime }/\mathbb{S}_{\nu },\mathcal{L}%
)$ of the Kontsevich compactification of the moduli spaces of Riemann
surfaces with coefficients in the local system $\mathcal{L}=Det(P_{\Sigma })$%
, where $P_{\Sigma }$ is the set of marked points. For any solution to the
quantum master equation (\ref{qms}) in $\widetilde{F}$, with $%
S_{0,1,1}=S_{0,1,2}=0$, the partition function on stable ribbon graphs $%
c_{S}(G)$ defines the characterestic homology classes in $H_{\ast }(%
\overline{\mathcal{M}}_{\gamma ,\nu }^{\prime }/\mathbb{S}_{\nu },\mathcal{L}%
)$.\vspace{-0.1in}
\end{proposition}

\section{Supersymmetric matrix integrals.}

\subsection{Odd inner product.}

Let $S(a_{i})$ be a solution to the quantum master equation (\ref{qms}) from 
$\widetilde{F}$, with $S_{0,1,2}=S_{0,1,1}=0$. Consider the vector space 
\begin{equation*}
M=Hom(V,End(U))
\end{equation*}%
where $\dim U=(d|d)$. The supertrace functional on $End(U)$ gives a natural
extension to $M$ of the odd symmetric inner product on $Hom(V,k)$ dual to $B$%
. Let us extend $S$ to a function $S_{gl}$ on $M$ so that each cyclically
symmetric tensor goes to the supertrace of the product of the corresponding
matrices from $End_{k}(U)$ 
\begin{equation*}
(a_{i_{1}},\ldots ,a_{i_{k}})^{c}\rightarrow tr(X_{i_{1}}\cdot \ldots \cdot
X_{i_{k}})
\end{equation*}%
and the product of cyclic words goes to the product of traces. The
commutator $I_{\Lambda }=[\Lambda ,\cdot ]$, for $\Lambda \in End(U)^{odd}$,
is an odd differentiation of $End(U)$. Notice that $I_{\Lambda }^{2}\neq 0$
for generic $\Lambda $. For such $\Lambda $ there always exists an operator $%
I_{\Lambda }^{-1}$ of regularized inverse: $[I_{\Lambda },I_{\Lambda
}^{-1}]=1$, preserving the supertrace functional. A choice of nilpotent $%
I_{\Lambda }^{-1}$is in one-to-one correspondence with $I_{\Lambda }^{2}-$%
invariant lagrangian subspace in $\Pi End(U)$, corresponding to $L=\{x\in
M|I_{\Lambda }^{-1}(x)=0\}$. Let $\Lambda =\left( 
\begin{array}{cc}
0 & Id \\ 
\Lambda _{01} & 0%
\end{array}%
\right) $, $\Lambda _{01}=diag(\lambda _{1},\ldots ,\lambda _{d})$, I take $%
I_{\Lambda }^{-1}\left( 
\begin{array}{cc}
X_{00} & X_{10} \\ 
X_{01} & X_{11}%
\end{array}%
\right) =\left( 
\begin{array}{cc}
M_{\lambda }X_{01} & M_{\lambda }(X_{00}+X_{11}) \\ 
0 & -M_{\lambda }X_{01}%
\end{array}%
\right) $ where $M_{\lambda }E_{i}^{j}=(\lambda _{i}+\lambda
_{j})^{-1}E_{i}^{j}$.

\begin{theorem}
\label{mxintodd}By the standard Feynman rules, the asymptotic expansion, at $%
\Lambda ^{-1}\rightarrow 0$, is given by the following sum over oriented\
stable ribbon graphs: 
\begin{equation*}
\log \frac{\int_{L}\exp \frac{1}{\hbar }\left( -\frac{1}{2}tr\circ
B^{-1}([\Lambda ,X],X)+S_{gl}(X)\right) dX}{\int_{L}\exp \frac{1}{\hbar }%
\left( -\frac{1}{2}tr\circ B^{-1}([\Lambda ,X],X)\right) dX}%
=const\sum_{G}\hbar ^{-\chi _{G}}c_{S}(G)c_{\Lambda }(G)
\end{equation*}%
where $\chi _{G}$ the euler characteristic of the corresponding surface, $%
c_{S}(G)$ is our partition function associated with the solution $S$, $%
c_{\Lambda }(G)$ is the partition function associated with $%
(End(U),tr,I_{\Lambda }^{-1})$ and constructed using the propagator $%
tr(I_{\Lambda }^{-1}\cdot ,\cdot )$. $c_{\Lambda }(G)$ defines the
cohomology class in $H^{\ast }(\overline{\mathcal{M}}_{\gamma ,\nu }^{\prime
}/\mathbb{S}_{\nu },\mathcal{L})$.
\end{theorem}

This follows from the standard rules (\cite{del}) of the Feynman
diagrammatics (compare with the formula (0.1) from \cite{GK}). In particular
the combinatorics of the terms in $S_{gl}(X)$ matches the data associated
with vertices in the complex $\mathcal{F}\widetilde{\mathbb{S}}$, i.e. the
symmetric product of cyclic permutations and the integer number. The
construction of $c_{\Lambda }(G)$ is studied in more details in (\cite[B3]%
{B3}). 

\subsection{Even inner product.}

Let now $S$ denotes a solution to quantum master equation in the
Batalin-Vilkovisky algebra (\ref{Feven}), with $S_{0,1,2}=S_{0,1,1}=0$. In
this case my basic matrix algebra is the general \emph{queer} superalgebra $%
q(U)$ (\cite{BL}) with its \emph{odd} trace. The associative superalgebra $%
q(U)$ is the subalgebra of $End(U\oplus \Pi U)$, 
\begin{equation*}
q(U)=\{X\in End(U\oplus \Pi U)|[X,\pi ]=0\}
\end{equation*}%
where $U$ is a purely even vector space and $\pi :$ $U\rightleftarrows \Pi U$
, $\pi ^{2}=1$ is the odd operator changing the parity. As a vector space $%
q(U)=End(U)\oplus \Pi End(U)$. The odd trace on $q(U)$ is defined as $otr(X)=%
\frac{1}{2}tr(\pi X)$. Let us extend $S$ to the function $S_{q}$ on $%
M=Hom(\Pi V,q(U))$, so that each cyclically symmetric tensor in $\Pi (\Pi
V^{\otimes j})^{\mathbb{Z}/j\mathbb{Z}}$ goes to the odd trace of the
product of the corresponding elements from $q(U)$ 
\begin{equation*}
(a_{i_{1}},\ldots ,a_{i_{j}})^{c}\rightarrow otr(X_{i_{1}}\cdot \ldots \cdot
X_{i_{j}})
\end{equation*}%
and the product of cyclic words goes to the product of the odd traces. Let
us denote by $otr\circ B^{-1}$ the odd extension to $M$, which is defined
using the \emph{odd} pairing $otr(XX^{\prime })$ on $q(U)$, of the even
symmetric inner product on $Hom(V,k)$. Let $\Lambda =\Pi diag(\lambda
_{1},\ldots ,\lambda _{N})$, $\Lambda \in \Pi End(U)$ be the odd element
from $q(U)$. The commutator $I_{\Lambda }=[\Lambda ,\cdot ]$ is an odd
differentiation of $q(U)$, and for generic $\lambda _{i}$, $I_{\Lambda
}^{2}\neq 0$, and $I_{\Lambda }$ is \emph{invertible} outside of the even
diagonal. Define the regularized inverse $I_{\Lambda }^{-1}$, $(I_{\Lambda
}^{-1})E_{i}^{j}=(\lambda _{i}+\lambda _{j})^{-1}\Pi E_{i}^{j}$, $%
(I_{\Lambda }^{-1})\Pi E_{i}^{j}=0$, $[I_{\Lambda },I_{\Lambda }^{-1}]=1$,
and let $L=\{x\in M|I_{\Lambda }^{-1}(x)=0\}$.

\begin{theorem}
By the standard Feynman rules, the asymptotic expansion, at $\Lambda
^{-1}\rightarrow 0$, is given by the following sum over oriented\ stable
ribbon graphs:%
\begin{equation*}
\log \frac{\int_{L}\exp \frac{1}{\hbar }\left( -\frac{1}{2}otr\circ
B^{-1}([\Lambda ,X],X)+S_{q}(X)\right) dX}{\int_{L}\exp \left( -\frac{1}{%
2\hbar }otr\circ B^{-1}([\Lambda ,X],X)\right) dX}=const\sum_{G}\hbar
^{-\chi _{G}}c_{S}(G)c_{\Lambda }(G),
\end{equation*}%
where $c_{S}(G)$ is the partition function from (\cite{B1}) and $c_{\Lambda
}(G)$ is the partition function, associated with $(q(U),otr,I_{\Lambda
}^{-1})$, and constructed using the propagator $otr^{dual}(I_{\Lambda
}^{-1}\cdot ,\cdot )$. $c_{\Lambda }(G)$ defines the cohomology class in $%
H^{\ast }(\overline{\mathcal{M}}_{\gamma ,\nu }^{\prime }/\mathbb{S}_{\nu })$%
, it is the generating function for the products of tautological classes.
\end{theorem}

\begin{theorem}
Let us consider the case of the one-dimensional vector space $V$ with even
symmetric inner product. The solutions in this case are arbitrary linear
combination of cyclic words $X^{3}$, $X^{5}$,\ldots , $X^{2n+1}$. Our matrix
integral in this case is a supersymmetric extension of the Kontsevich matrix
integral from (\cite{K}).\vspace{-0.05in}
\end{theorem}

\end{document}